\documentstyle[amssymb,12pt]{article}

\newtheorem{def-prop}{Definition-Proposition}
\newtheorem{thm}{Theorem}

\newtheorem{lemma}{Lemma}
\newtheorem{defn}{Definition}
\def\endproof{$\Box$}
\def\1{1}
\def\P{{\bf P}}

\def\C{{\Bbb C}}
\def\Z{{\Bbb Z}}
\def\A{{\Bbb A}}

\def\Q{{\Bbb Q}}

\def\co{{\cal O}}

\begin{document}
\title{Algebraic cuts}
\author{Dan Edidin\thanks{Partially supported by the NSF,
and the University of Missouri Research Board.}
\hspace{1mm} and
William Graham}
\date{}
\maketitle

\section{Introduction}

Let $X$ be a projective variety with a linearized action of
an algebraic group $G$.  If $X$ is smooth, and the ground field
is $\C$, then the geometric invariant theory quotient $X//G$ can be
identified with a quotient constructed using symplectic geometry, the
``reduced space'' $X_r$.  This result, due to Mumford, Guillemin
and Sternberg, connects geometric invariant theory and
symplectic geometry.

Suppose now that $G = T$ is a torus, which for simplicity we will take
to have dimension 1.  In \cite{Lerman}, Lerman introduced a construction
called symplectic cutting, which constructs a manifold $X_c$ related to
$X$, with a $T$-action, and embeds $X_r$ as a component of the fixed
point locus $X_c^T$.  The complement of $X_r$ in $X_c$ can be identified
with an open submanifold of $X$, so the other components of $X_c^T$ are
certain components of $X^T$.  The components of $X_c^T$ are linked by
the localization theorem in equivariant cohomology.  Thus, from
knowledge of $X^T$ one can (via the cut space) deduce results about
$X_r$.  For example, Lerman uses cutting to prove a residue formula due
to Kalkman, which is closely connected to the localization theorem of
Jeffrey-Kirwan-Witten (\cite{G-K}).

The purpose of this paper is to present an algebraic version of Lerman's
construction, called algebraic cutting, which is valid over arbitrary
ground fields and for possibly singular schemes.  This is useful in
studying $X_r$ from the point of view of algebraic geometry, using the
equivariant intersection theory developed in \cite{E-G} in place of
equivariant cohomology.  For example, Lerman's proof of Kalkman's
formula becomes valid for smooth schemes over arbitrary ground fields;
at the end of this paper we briefly give an adaptation of that proof.

Throughout this paper $T$ will denote a one-dimensional split torus
${\Bbb G}_m$, $X \subset \P(V)$ a quasiprojective variety with a
linearized $T$-action, and $X^s$ the set of stable points for this
action.  An element of the group $T$ will usually be denoted by $g$.
All Chow groups are assumed to have rational coefficients.  $A^*_T$
denotes the $T$-equivariant Chow ring of a point; it is isomorphic to
the polynomial ring $\Q[t]$.  The sign convention for $A^*_T$ is
determined as follows.  We model $BT$ by $(V - {0})/T = \P^N$, where $V$ is an
$N+1$ dimensional vector space on which $T$ acts with all weights $-1$.
The element $t$ represents the hyperplane class.  A representation of
$T$ is an equivariant vector bundle over a point, and thus has
equivariant Chern classes, coming from the induced bundle over $BT$.  A
1-dimensional representation with weight $a$ then has equivariant first
Chern class $at$.  The Chow ring $A^*_{T \times T}$ will be denoted
$\Q[t_1,t_2]$.

\section{Algebraic cuts}

\subsection{Cutting at zero}
Let $v_1 , \ldots, v_n$ be a basis of $V$ consisting of
weight vectors, of weights $a_1 , \ldots, a_n \in \Z$.  If $x \in \P(V)$
is the line spanned by $\sum x_i
v_i$, define (cf. \cite{B-P}) the weights of $x$
(denoted $\Pi(x)$) to be the set of $a_i$ such that $x_i$ is nonzero.
As an algebraic analogue of the fact that 0 is
a regular value of the moment map, we will assume that for all $x
\in X$, $0 \notin \Pi(x)$.
Thus, $X$ actually lies in $\P(W)$, where $W$ is the
subspace of $V$ spanned by weight vectors with nonzero weights.  By
replacing $V$ by $W$ if necessary, we will also assume that $0$ is not a
weight of $T$ on $V$.

By \cite{B-P} stability with respect to the line bundle $\co_{\P(V)}(1)$
can be expressed in terms of the weights:
$$X^s = \{x \in X | \Pi(x) \mbox{ contains both positive and negative
  weights} \}.$$ By assumption, $0$ is not a weight of $T$ on $V$, so
$X^{ss} = X^{s}$.  Hence if $X$ is projective then the ``reduced space''
$X_r := X_s/T$ is projective.  If $X$ is smooth and the ground field is
$\C$ then $X_r$ corresponds to the symplectic reduction at zero.

Define
$$
X_{>0} = \{x \in X | \Pi(x) \mbox{ contains a positive weight} \}.
$$
Let $T_{\Delta}$ and $A$ denote the diagonal and antidiagonal copies
of $T$ in $T \times T$.  Let $(X \times \A^1)^s$ denote the stable
points of $X \times \A^1$ with respect to the action of $A$ and the
line bundle ${\cal L} = \co_{\P(V)}(1) \boxtimes 1$.

\begin{defn} \label{d.1}
The algebraic cut of $X$ is the quotient $X_c = (X
\times \A^1)^s /A$.  This scheme has a natural $T$-action, coming from
the isomorphism $T \rightarrow (T \times T)/A$, $g \mapsto (g,1) \mbox{
  mod }A$.
\end{defn}

The main results concerning this space are two theorems,
which are algebraic versions of Lerman's results.  The first
says that $X_c$ is projective (if $X$ is) and is assembled out of the pieces
$X_{>0}$ and $X_r$.

\begin{thm} \label{t.cuts1}
$(a)$  If $X$ is projective then $X_c$ is  projective.

$(b)$ There is a $T$-equivariant open embedding $i_c: X_{>0}
\hookrightarrow X_c$ and a $T$-equivariant closed embedding $j_c: X_r
\hookrightarrow X_c$ such that $X_c$ is a disjoint union of $X_{>0}$
and $X_r$.  (Here $T$ acts trivially on $X_r$.)

$(c)$ A neighborhood of $X_r$ in $X_c$ is isomorphic to $X^s \times^T
\A^1 \rightarrow X_r$.

\end{thm}

Proof: To prove (a), we will use the fact that
if $Y$ is a projective variety with a linearized $G$-action, then
$Y^{ss}/G$ is projective (\cite{GIT}).  Since $X \times \A^1$ is not
projective we will construct an $A$-equivariant embedding into a
projective variety $Y$ with a linearized $A$-action and show that
$Y^{ss}$ coincides with $(X \times
\A^1)^s$.  We do this as
follows.  For any integer $j$ let $\A^2_{j}$ denote $\A^2$ with weights
$0$ and $j$ for the action of $T$. Likewise, set $\P^1_{j} = \P(A^2_j)$
with the induced $T$-action.  The embedding $\A^1 \hookrightarrow
\P^1_{1}$ given by $z \mapsto [1:z]$ is $T$-equivariant and induces a $T
\times T$-invariant embedding $p: X \times \A^1 \hookrightarrow X \times
\P^1_{1}$.

Let $a_0, \ldots , a_n$ be the weights of $T$ on $X$.
Fix
an integer $N$ which is greater than any $a_i$.  There is a product
action of $T \times T$ on $V \otimes \A^2_{N}$. The weights
$T_{\Delta}$ on $V \otimes \A^2_{N}$ are $a_0,\ldots,a_n,a_0+N,\ldots,
a_n+N$, and the weights of $A$ are $a_0,\ldots,a_n,a_0-N,\ldots,
a_n-N$. Since $N$ is greater than each $a_i$,
$0$ is not a weight of $A$ on $V \otimes \A^2_{N}$. Thus,
by \cite{B-P} all semi-stable points
(with respect to $A$) are stable.

Consider the $T \times T$-equivariant finite morphism
$\psi : X \times \P^1_{1} \rightarrow \P(V \otimes \A^2_{N})$
given by
$$
([x_0: \cdots :x_n], [z_0:z_1]) \mapsto [x_0 z_0^N: \cdots :x_n z_0^N :
x_0 z_1^N: \cdots x_n z_1^N].
$$
Let ${\cal M} = \psi^* \co (1)$.

Let $\P(V \otimes \A^2_{N})^s$ denote the stable points for the action
of $A$ with respect to $\co (1)$ and $(X \times \P^1_{1})^s$ denote
the stable points for the action of $A$ with respect to ${\cal M}$.
Since
$\psi$ is finite and $X \times \P^1_{1}$ is proper, by \cite[Theorem
1.19]{GIT}, $(X \times \P^1_{1})^s = \psi^{-1}(\P(V \otimes
\A^2_{N})^s)$. Since 0 is not a weight, all semi-stable points
are stable.

The stable points in $\P(V \otimes \A^2_{N})$ are described in
\cite{B-P}: if we let $\Pi_A(y)$ denote the $A$-weights of $y \in \P(V
\otimes \A^2_{N})^s$, then $y$ is stable if and only if $\Pi_A(y)$
contains both positive and negative weights.  Hence
the unstable points $(x,z)$ of $X \times \P^1_1$ with respect to ${\cal M}$ are
of 3 types:

$(i)$ $z=[0:1]$, $x$ arbitrary;

$(ii)$ $z=[1:0]$, $\Pi(x)$ all positive or all negative;

$(iii)$ $z \neq [0:1],[1:0]$, $\Pi(x)$ all negative.\\ In particular,
the stable points of $X \times \P^1_1$ lie in the open set $X \times
\A^1 \subset X \times \P^1$.  Restricting to points $(x,w)$ in the open
set $X \times \A^1$ we see that there are two types of unstable points:

$(ii')$ $w=0$, $\Pi(x)$ all positive or all negative;

$(iii')$ $w\neq 0$, $\Pi(x)$ all negative.

\medskip

The restriction of ${\cal M}$ to $X \times \A^1$ is ${\cal L}$.  As in
Definition \ref{d.1}, let $(X
\times \A^1)^s$ denote the set of stable points
relative to to $T$ and ${\cal L}$.  By \cite[Prop. 1.18]{GIT} the inclusion
$p:X \times \A^1 \subset X \times \P^1_{1}$ induces a reverse inclusion
$$
(X \times \A^1)^s \supset p^{-1}(X \times \P^1_{1})^s = (X \times
\P^1_{1})^s \cap (X \times \A^1) .
$$
However, since $(X \times \P^1_1)^s \subset X \times \A^1$
this implies
$$
(X \times \A^1)^s \supset (X \times \P^1_{1})^s.
$$
To prove that $X_c = (X \times \A^1)^s/A$ is projective
it suffices to show that this inclusion is an equality, since the
quotient $(X \times \P^1_{1})^s /T$ is projective.
In other words, we must show that
all points of type (ii') and (iii') are ${\cal M}$-unstable.
We will do this the old fashioned way, and show that any
$A$-invariant element of $H^0(X \times \A^1, {\cal L}^d)$ vanishes on points
of type (ii') and (iii').
In doing so, we may assume $d >>0$.

For $d>>0$, the restriction map
$$
H^0(\P(V),{\cal O}(d)) \cong S^d(V^*) \rightarrow H^0(X,{\cal O}(d)),
$$
is surjective (here $V^*$ is the dual vector space to $V$).  If $\xi \in
V^*$ is a positive (resp. negative) weight vector and if $x \in X$
has all positive (resp. negative) weights,
then $\xi$ (viewed as a section of ${\cal O}(1)$)
vanishes at $x$.

Now if $s \in S^d(V^*)$ is a $T$-invariant section, then
it is a sum of monomials in weight vectors, with each monomial containing both
positive and negative weight vectors. Thus, if $s \in S^d(V^*)^T$ and
$\Pi(x)$ consists of all positive or all negative weights, then $s(x)
= 0$. Similarly, if $s \in S^d(V^*)$ is a sum of monomials and each
monomial contains a negative weight vector (in particular if the weight
of $s$ is negative), then $s(x) = 0$ if $\Pi(x)$ has only negative weights.

Now if $s \in H^0(X,{\cal O}(d))$ has $T$-weight $k$, and
$f_l(w) = w^l \in {\cal O}(\A^1)$, then $s \otimes f_l$ has weight
$k+l$ for the anti-diagonal torus $A$ (note that $f_l$ has $T$-weight
$-l$).
Thus, a general $A$-invariant
section $\tau$ of $H^0(X \times \A^1, {\cal L}^d)$ is of the
form $\tau = \sum \tau_lf_l$, where $\tau_l \in H^0(X,{\cal O}(d))$ has
$T$-weight $-l$.

We want to show that if $\tau$ is $A$-invariant then
$\tau$ vanishes on $(x,w)$ of type $(ii')$ and $(iii')$.  On type
$(ii')$, we have
$$
\tau(x,0) = \sum \tau_l(x) f_l(0) = \tau_0(x).
$$
But $\tau_0$ is $T$-invariant, and $(x,0)$ of type
$(ii')$ means that $\Pi(x)$ consists of all positive or all negative
weights, so by the above remarks, $\tau_0(x) =0$.  If $(x,w)$ is of
type $(iii')$, then $\Pi(x)$ consists of all negative
weights; since the weight of $\tau_i = -i \leq 0$, we have $\tau_i(x)
=0$ for all $i$.  Hence $\tau$ vanishes on $(x,w)$ of type $(ii')$ and
$(iii')$, as desired. This proves (a).

We now prove (b) and (c).
The embedding $i_c$ is defined to be the quotient by $A$ of the open embedding
$$
X_{>0} \times T \hookrightarrow (X \times \A^1)^s;
$$
The embedding $j_c$ is defined to be the quotient by $A$ of the
embedding
$$
J: X^s \times \{ 0 \} \hookrightarrow (X \times \A^1)^s;
$$
$J$ is closed since, from the description of $(X \times \A^1)^s$
given in the proof of (a), we have $(X \times \A^1)^s \cap X \times \{
0 \} = X^s \times \{ 0 \}$.
The embedding $i_c$ (resp. $j_c$) is open (resp. closed) since it is the
quotient of an open (resp. closed) embedding.  This proves (b).

A neighborhood of $X^s \times \{ 0 \}$ in $(X \times \A^1)^s$ is given
by $X^s \times \A^1$.  It follows that a neighborhood of ${X_r}$ in $X_c$
is $T$-equivariantly isomorphic to the fibration
$$
(X^s \times \A^1)/A \rightarrow (X^s \times \{ 0 \})/A.
$$
This is isomorphic to the fibration
$$
X^s \times^T \A^1 \rightarrow X_r,
$$
proving (c). \endproof

\medskip

{\bf Remark.} In general $T$ will act with finite stabilizers on $X^s$.
If $T$ acts freely on $X^s$ then $X_r$ is regularly embedded with normal
bundle $X^s \times^T \A^1 \rightarrow X_r$.  We will see in Theorem
\ref{t.cuts2} that even if the action of $T$ is not free there is an
equivariant pullback and self-intersection formula.

\subsection{Cutting at any point} \label{cutany}

Instead of considering stable points with respect to the invertible sheaf
$\co_{\P(V)}(1)$ we can consider stable points with respect to other
linearizations.
Let $q = a/n$ be a rational number which is not a weight of $\Pi(x)$
for any $x \in X$ (we no longer require that $0$ is not a weight).  Then
we define
$$
X^s(q) = \{x \in X | \Pi(x) \mbox{ contains both weights greater than
  $q$ and weights less than $q$} \}.
$$
Define a $T$-action on $S^n V$ by $g \cdot v^n = g^{-a} (g \cdot
v)^n$.  Embed $\P(V)$ in $\P(S^n V)$ by the Veronese embedding, and let
${\cal M}(q)$ denote the pullback to $X$ of $\co_{\P(S^n V)}(1)$.  Then
(\cite[1.2]{B-P}) $X^s(q)$ is the set of stable points of $X$ with
respect to ${\cal M}(q)$.  The quotient $X_r(q) = X^s(q)/T$ corresponds
(in the symplectic picture) to the reduction of $X$ at $q$.

Likewise we can define (in the obvious fashion) the open subscheme
$X_{>q}$ of $X$.  We let ${\cal L}(q)$ denote the line bundle ${\cal
  M}(q)\boxtimes 1$ on $X \times \A^1$, and define $X_c(q) = (X
\times \A^1)^s /A$, where now stability is defined with respect to
${\cal L}(q)$.   The analogue of Theorem \ref{t.cuts1} holds for $X_c(q)$.
This can be deduced from Theorem \ref{t.cuts1} by embedding $X$ into
$\P(S^n V)$ as above.

\section{Equivariant Intersection theory and algebraic cuts}
\subsection{Review of Equivariant Chow groups}

In this section we recall the definition and some of the
basic properties of equivariant
Chow groups \cite{E-G}. For simplicity we will
assume that all schemes are quasiprojective and
that any group actions are linearized.

Let $G$ be a $g$-dimensional group, $X$ an $n$-dimensional scheme
and $V$ a representation of $G$ of dimension $l$. Assume that there
is an open set $U \subset V$ such that a principal bundle
quotient $U \rightarrow U/G$ exist, and that $V-U$ has codimension
more than $i$. Let $X_G = (X \times U)/G$.
\begin{defn}
Set $A_i^G(X) = A_{i+l-g}(X_G)$, where $A_*$ is the usual Chow group.
This definition is independent of the choice of $V$ and $U$ as long as
$V-U$ has sufficiently high codimension.
\end{defn}

{\bf Remark:} Because $X \times U \rightarrow X \times^G U$ is
a principal $G$-bundle, cycles on $X \times^G U$ exactly
correspond to $G$-invariant cycles on $X \times U$. Since
we only consider cycles of codimension smaller
than the dimension of $X \times (V-U)$, we may in fact
view these as $G$-invariant cycles on $X \times V$. Thus
every class in $A_i^G(X)$ is represented by a cycle in
$Z_{i+l}(X \times V)^G$, where $Z_*(X \times V)^G$ indicates
the group of cycles generated by invariant subvarieties.
Conversely, any cycle in $Z_{i+l}(X \times V)^G$ determines
an equivariant class in $A_i^G(X)$.

\medskip

The properties of equivariant intersection theory include the following.

(1) Functoriality for equivariant maps:
proper pushforward, flat pullback, l.c.i pullback, etc.

(2) Chern classes of equivariant bundles operate on  equivariant
Chow groups.

(3) If $X$ is smooth, then $\oplus A_*^G(X)$ has a ring structure. (This
follows from (1), since the diagonal $X \hookrightarrow X \times X$
is an equivariant regular embedding when $X$ is smooth.)

(4) There is a localization theorem for torus actions.

\begin{lemma} \label{l.normal}
  Let $G$ be a group and $X$ a quasi-projective $G$-scheme. Suppose that
  $X$ is the set of stable points for the linear action of a normal
  (reductive) subgroup $H \subset G$. Then $A_i^G(X) =
  A_{i-h}^{G/H}(X/H)$, where $h = \mbox{dim }H$.
\end{lemma}

Proof. Let $V$ be a representation of $G/H$, and let $U$ be an open
set on which $G/H$ acts freely. Then $G$ acts on $X \times U$ with finite
stabilizers. Thus, by \cite[Theorem ?]{E-G}  $A_*^G(X \times U) =
A_*( (X \times U)/G)$. Now $(X \times U)/G = G/H(X \times U)/H$.
Since $H$ acts trvially on $U$, $(X \times U)/H = X/H \times U$.
Thus $(X \times U)/G = X/H \times^{G/H} U$.

Hence if $V - U$ has sufficiently high codimension then
$A_*(X \times U) = A_*^{G/H}(X/H)$. Thus, the lemma follows from the
homotopy property of equivariant Chow groups.
\endproof

\begin{lemma}  Suppose $X$ is a $G$-scheme, $Y$ is an $H$-scheme, and
$Q \subset X \times Y$ is $G \times H$-invariant and flat over $X$.
Let $\pi: Q \rightarrow X$ be the projection.
Then there is a pullback
$$
\pi^* : A^G_*(X) \rightarrow A_*^{G \times H}(Q)
$$
\end{lemma}
Proof. Let $V_G$ be a representation of $G$. Let $p_{Q}: Q \times V_G
\rightarrow X \times V_G$
be the projection onto the first and third factors.
If $Z \subset (X \times V_G)$ set $\pi^*_Q[Z] = [\pi_{Q}^{-1}Z]$.
This is an invariant cycle in $Q \times V_G$ and so defines an
element of $A_*^{G \times H}(Q)$. The usual arguments of equivariant
theory show that this is well defined.
\endproof

\begin{thm} \label{t.cuts2}
Let $X$ be a quasiprojective scheme with a linearized
  $T$-action and cut scheme $X_c = X_c(q)$ for some $q \in \Q$. Then:

(a) There is a pullback of rational equivariant Chow groups
$$j_c^*: A_*^T(X_c) \rightarrow A_*^T(X_r)$$
such that
$j_{c*}  j_c^*$ is multiplication by
$t$ under the map $A^1_T \rightarrow  A^1_T(X^s) \cong A^1(X_r)$.

$(b)$ There is a map $s: A_*^T(X) \rightarrow A_*^T(X_c)$ such that
$i^*_c \circ s = i^*$.  Here $i$ and $i_c$ denote the embeddings of
$X_{>0}$ into $X$ and $X_c$, respectively.

$(c)$ If $F \subset X_{>0}^T$, let $i_F$ and $i_{F,c}$ denote the
embeddings of $F$ into $X$ and $X_c$ respectively.  Then $i_{F,c}^*
\circ s = i_F^*$.

$(d)$ As maps $A_*^T(X) \rightarrow A_*(X_r)$, we have $ {\cal F} \circ
j_c^* \circ s = r^*$, where ${\cal F}$ denotes the forgetful map from
equivariant Chow groups to ordinary Chow groups, and $r^*$ is the
restriction $A_*^T(X) \rightarrow A_*^T(X_s) \cong A_*(X_r)$.
\end{thm}

{\bf Remark.}  Because $T$ acts trivially on $X_r$, we have $A^*_T(X_r)
\cong A_*(X_r) \otimes_{\Q} A^T_*$.  Thus, if $\alpha \in A_*^T(X)$, we
can write $j_c^* \circ s(\alpha) = \sum \beta_i t^i$.  The content of
(d) is that $\beta_0 = r(\alpha)$.

\medskip

Proof:
To simplify notation we assume $q = 0$.
By the discussion of Section \ref{cutany} this implies the theorem for
$q \neq 0$ as well.
(a) The $T \times T$ equivariant regular embedding $J: X^s \times
\{0\} \hookrightarrow (X \times \A^1)^s$ gives a pullback $J^*$ on $T
\times T$-equivariant Chow groups.  As noted in the proof of Theorem
\ref{t.cuts1}, the map $j_c$ is the quotient of the map $J$ by $A$.  We
have an isomorphism $T \rightarrow (T \times T) / A$, $g \mapsto
(g,1) = (1,g)$ mod $A$.  By Lemma \ref{l.normal} we can identify
$$
A_*^T(X_c) = A_*^{T \times T}((X \times \A^1)^s) \mbox{   ,   }
A_*^T(X_r) = A_*^{T \times T}(X^s \times \{0\}).
$$
Let $j_c^*$ correspond
to $J^*$ under this identification.  The $T \times T$ equivariant normal
bundle to $X^s \times \{0\}$ in $(X \times \A^1)^s$ is the trivial line
bundle, where $T \times T$ acts with weight $(0,1)$.  Thus $J_*J^*$ is
multiplication by $t_2 \in A^*_{T \times T}$.  With our identifications,
multiplication by either $t_1$ or $t_2$ in $A^*_{T \times T}=
\Q[t_1,t_2]$ corresponds to multiplication by $t \in A^*_T = \Q[t]$.
Hence $j_{c *}j_c^*$ is multiplication by $t$. This proves
(a).

(b) We have inclusions
$$
X_{>0} \times T \stackrel{I} \rightarrow (X \times \A^1)^s \stackrel{k}
\rightarrow X \times \A^1 .
$$
The map $i_c$ is the quotient of $I$ by $A$.  Consider the following
commutative diagram:
$$
\begin{array}{cccc}
A_*^T(X) \stackrel{\pi^*} \rightarrow A_*^{T \times T}(X \times \A^1)
\stackrel{k^*} \rightarrow & A_*^{T \times T}((X \times \A^1)^s) &
\stackrel{I^*} \rightarrow & A_*^{T \times T} (X_{>0} \times T) \\
  & \small{f} \downarrow &                   & \small{g} \downarrow \\
  & A_*^T(X_c)   &\stackrel{i_c^*} \rightarrow  & A_*^T( X_{>0} )
\end{array}
$$
The vertical maps are from Lemma \ref{l.normal}; here $G = T \times T$,
$H = A$, and $G/H = (T \times T)/A \cong T$.  We define
$$
s: A_*^T(X) \rightarrow A_*^T(X_c)
$$
to be $s = f k^* \pi^*$.  If $Z$ is a subvariety of $X$, then

$$
i_c^* \circ s ([Z]) = g I^* k^* \pi^*[Z] = g[(Z \cap X_{>0}) \times T] .
$$
Since $A$ acts freely on $X_{>0} \times T$,
$$
g[(Z \cap X_{>0}) \times T] = [Z \cap X_{>0}] = i^*[Z] ,
$$
proving (b).  Part (c) follows from (b).

The proof of (d) is similar to (b).  \endproof

\subsection{Kalkman's residue formula}
Following \cite{Lerman}, we use algebraic cutting to
prove the residue formula of Kalkman.

\begin{thm} \label{t.kalkman}
Let $X \subset \P(V)$ be a smooth $n$-dimensional projective variety
with
a linearized $T$-action. Fix a rational $q$ which is not a weight of
$x \in X$.
Let $r^*: A^*_T(X) \rightarrow A^*(X_r)$ be as defined above, and let
$\alpha \in A^{n-1}_T(X)$.  Then
$$
\mbox{deg } (r^*(\alpha) \cap [X_r]) = - \sum \mbox{deg }( (\mbox{Res}_{t=0}
\frac{i^*_F \alpha}{c^T_{d_F}(N_F X)}) \cap [F]_T) ,
$$
where the sum is over the components $F$ of $X^T_{>q}$, and $d_F$ is
the codimension of $F$.
\end{thm}

Proof:
To simplify notation assume that $q = 0$ in
the proof.
Recall that $X^T_c = X^T_{>0} \cup X_r$.  Each $F \subset
X^T_{>0}$ is regularly embedded in $X_c$, so the self-intersection
formula applies.  Although $X_r$ need not be regularly embedded in
$X_c$, a self-intersection formula still holds (Theorem
\ref{t.cuts2}(a)).  Thus the localization theorem for equivariant Chow
groups \cite{E-G2} can be applied just as if $X_r$ were regularly
embedded, so
\begin{equation} \label{e.kalkman1}
s(\alpha) \cap [X_c]_T = \sum_{F \subset X_{>0}}
i_{F,c*}( \frac{i^*_{F,c} s(\alpha)}{c^T_{d_F}(N_F X)} \cap [F]_T )+
j_{c*} ( \frac{j_c^*s(\alpha)}{t} \cap [X_r]_T)  .
\end{equation}
This equality holds in $A_*^T(X_c) \otimes_{A_T^*} {\cal Q}$, where
${\cal Q} = \Q[t,t^{-1}]$.  Applying
Theorem \ref{t.cuts2}(c) we can rewrite this as
\begin{equation} \label{e.kalkman2}
s(\alpha) \cap [X_c]_T = \sum_{F \subset X^T_{>0}}
i_{F,c*} ( \frac{i_F^*\alpha}{c^T_{d_F}(N_F X)}  \cap [F]_T) + j_{c*}
( \frac{j_c^* s(\alpha)} {t} \cap [X_r]_T).
\end{equation}
Let $\pi_Y$ denote the projection of a proper scheme $Y$ to a point and
$\pi_{Y*}$ the induced map of equivariant Chow groups.  If we apply
$\pi_{X_c*}$ to both sides of (\ref{e.kalkman2}),
the left side
is zero since $\alpha$ has degree $n-1$ and the dimension of
$X$ is $n$.  Therefore
\begin{equation} \label{e.kalkman3}
\sum_{F \subset X^T_{>0}}
\pi_{F*}( \frac{i_F^*\alpha}{c^T_{d_F}(N_F X)}  \cap [F]_T)=
- \pi_{X_r*}( \frac{j_c^* s(\alpha)}{t} \cap [X_r]_T).
\end{equation}
Now, $T$ acts trivially on $F$, so $A^*_T(F) = A^*F \otimes_{\Q} A^*_T$
and $[F]_T = [F] \otimes 1$.  The class
$$
\frac{i^*_F \alpha}{c^T_{d_F}(N_F X)} \in A^*_T(F) \otimes_{ {\cal Q}}
A^*_T = A^*F \otimes_{\Q} {\cal Q}
$$
has degree $n-1-d_F = \mbox{dim }F - 1$, so we can write it as $\sum
\beta_i t^i$, where $\beta_i \in A^*F$ has degree $\mbox{dim }F - 1 -
i$.  Thus
$$
\frac{i^*_F \alpha}{c^T_{d_F}(N_F X)} \cap [F]_T = (\sum \beta_i \cap
[F]) t^i.
$$
When we apply $\pi_{F*}$ to this the only term that survives is the term
with $t^{-1}$.  A similar argument applies to the right side of
(\ref{e.kalkman3}).  We can therefore rewrite (\ref{e.kalkman3}) as
$$
 \sum_{F \subset X_{>0}^T} \mbox{deg }( (\mbox{Res}_{t=0}
\frac{i^*_F \alpha}{c^T_{d_F}(N_F X)}) \cap [F])
= \mbox{deg } \mbox{Res}_{t=0} \frac{j^* s(\alpha) \cap [X_r]_T}{t}.
$$ By the remark after Theorem \ref{t.cuts2}, the right side of this
equation is $r^*(\alpha) \cap [X_r]$.
\endproof

\subsection{Characteristic numbers of quotients}  Kalkman's residue
formula has a number of applications; for instance it can be used to
prove the localization theorem of Jeffrey-Kirwan-Witten \cite{G-K}.  Another
use of the residue formula is to compute some characteristic
numbers of quotients of a smooth variety $X$ by a torus.  We present some
formulas when $T$ acts freely on $X^s(q)$ and the fixed points
of $X_{> q}$ are isolated.

Let $Y = X^s/T$ and let $\pi: X^s \rightarrow Y$ be the quotient map.
Since $T$ acts freely, $Y$ is smooth of dimension
$n-1$ and
we will use the residue formula to compute the Euler characteristic
$\chi(Y) = \mbox{deg }( c_{n-1}(T_Y))$.
Since the fiber of $\pi$ is $T$, $T_{X^s/Y}$ is the trivial bundle
of rank $1$. There is an exact sequence of equivariant bundles on
$X^s$
$$ 0 \rightarrow {\bf 1} \rightarrow T_{X^s} \rightarrow \pi^*{T_{Y}}
\rightarrow 0.$$
Hence $r(c_i(T_{X^s})) = c_i(T_{Y})$. Since $T_{X^s}$ is just
the restriction of $T_X$ to $X^s$, we can apply the residue formula.
If $p \in X_{>q}^T$ is a fixed point let
$\alpha_{1}(p), \ldots , \alpha_{n}(p)$ be the weights for the $T$-action
on $T_{p}X$. Then
$$c^T_{n-1}(T_{X}|_{p}) = t^{n-1} \sum_{i = 1}^{n} \alpha_1(p) \alpha_2(p)
\ldots \widehat{\alpha_{i}(p)} \ldots \alpha_n(p)$$ and
$c^T_{n}(N_{p}X) = t^n\prod_{i=1}^n \alpha_i(p)$.
Thus applying the residue formula we have
$$\chi(Y) = - \sum_{p \in X_{>q}^T} \sum_{i = 1}^{n} \frac{1}{\alpha_i(p)}.$$

In a similar manner we can calculate the Todd genus $\chi({\cal O}_Y)$.
By Hirzebruch-Riemann-Roch,
$\chi({\cal O}_Y)= \mbox{deg }(Td(T_{Y})).$
Using
the residue formula to calculate $\mbox{deg }(Td(T_{Y}))$
we obtain
$$\chi({\cal O}_Y) = - \sum_{p \in X_{>q}^T} \mbox{res}_{t =
  0}\left(\frac{1}{ \prod_{i =1}^{n}(1 - e^{-t \alpha_i(p) })}\right) =
- \sum_{p \in X_{>q}^T} \left(\sum_{1 \leq i \neq l \leq n}
\frac{\alpha_{l}(p)} {2\alpha_i(p)}\right)$$.

\medskip

\noindent { Deparment of Mathematics, University of Missouri, Columbia MO
65201}\\

\noindent { Insitute for Advanced Study, Princeton NJ 08540}
\end{document}